\documentclass[final]{jcmlatex}

\usepackage{amssymb}
\usepackage{amsmath}
\usepackage{amsfonts}
\usepackage{graphicx}

\newcommand{\xvec}{\mathbf x}

\newcommand{\dnx}{\text{d}^n\xvec}
\newcommand{\mathlim}[2]{\underset{#1 \rightarrow #2}{\text{lim}}}
\begin{document}

\title{A global root-finding method for high dimensional problems}
\author{Fabrizio Castellano
               \thanks{Dipartimento di Fisica, Politecnico di Torino, Corso
Duca degli Abruzzi 24, 10129 Torino, Italy\\
Email: fabrizio.castellano@polito.it}
    }

\maketitle

\begin{abstract}
A method to solve the problem $f(x) = 0$ efficiently on any $n$-dimensional domain $\Omega$ under very broad hypoteses is proposed. The position of the root of $f$, assumed unique, is found by computing the center of mass of an $\Omega$-shaped object having a singular mass density. It is shown that although the mass of the object is infinite, the position of its center of mass can be computed exactly and corresponds to the solution of the problem. The exact analytical result is implemented numerically by means of an adaptive Monte Carlo sampling technique which provides an exponential rate of convergence. The method can be extended to functions with multiple roots, providing an efficient automated root finding algorithm.
\end{abstract}

\section{Introduction}
Needless to say, the problem of finding the zeroes of a function is indeed one of the central issues in every branch of science, since many others can be reduced to it. 
Scalar equations of the form $f(\xvec) = g(\xvec)$, where $\xvec$ is an $n$-dimensional vector, can be rewritten as $f(\xvec) - g(\xvec) = 0$.
Nonlinear systems of equations $f_\alpha(\xvec) = 0$ can be rewritten as $F(\xvec) = 0$, where $F$ is a vector valued function, and then reduced to a scalar equation by taking $|F(\xvec)| = 0$.

When talking about finding the root of a function the dimension of the space plays a crucial role. For functions in one variable well estabilished and efficient algorithms exist, like the bisection, secant, Newton-Raphson methods and many others. These methods can either be higly efficient in terms of rate of convergence but may fail to converge in certain cases or can ensure convergence at the cost of a slower rate. The key feature in one dimension is that it is always possible to bracket the solution between two points, use such points as an estimate for the root and then improve above them. This is what allows to develop convergent algorithms.
In two or more dimensions the problem is completely different because the boundary of any domain containing the solution is made up of an infinity of points rather than only two.
Thus multidimensional root-finding problems are usually deterministically approached by methods involving local quantities, such as a guess of the solution and the derivatives of the function around it. These are extensions of the Newton-Raphson method\cite{basics} and may fail to converge if the function is not sufficiently well behaved in the domain of interest: rapidly oscillating or not-differentiable functions pose significant problems in the application of such methods.

As an alternative to deterministic methods, stochastic root finding algorithms provide an efficient way to cope with bad bahaved functions, and are the only available tool when the value of the function itself cannot be known exactly without a high computational cost\cite{SRFP}.
Another promising alternative to root finding and other optimization problems is represented by genetic algorithms and evolutionary programming \cite{GA1, GA2}, in which the ``survival of the fittest'' paradigm is used to select the solution among a population of initial candidates which are progressively mutated and discarded basing on a suitable fitness function.

One of the key problems affecting all the above-mentioned methods is that at least one initial guess of the solution is required, from which the method starts. Such initial guess is what ultimately determines if the method will converge to the correct solution or will get trapped inside a local minumum or other bad behaved regions.
This weakness is due to the intrinsic locality of such methods, which try to improve on the available solution basing on informations regarding only a neighborhood of the current estimate.

In the present paper an exact analytical result that allows to directly compute the position of the root of a function is presented. Basing on this result a numerical method to find the solution is developed, which is intrinsically global and does not need an initial guess of the solution.

\section{The center of mass property}
Let us consider a scalar function $f$ defined over an $n$-dimensional domain $\Omega$ having a unique zero at $\xvec_0 \in \Omega$. We will further assume that there exists a simply-connected open set $B^\varepsilon_{\xvec_0} \subset \Omega$ ,of arbitrarily small characteristic size $\varepsilon$, such that $\xvec_0 \in \partial B^\varepsilon_{\xvec_0}$ and over which $f$ is differentiable. 
In practice we require $f$ and $\Omega$ to be well behaved just over an arbitrarily small domain containing the root on its boundary.
Given $f$ and $\Omega$ satisfying the above conditions let us consider the family of functions
\begin{equation}
\label{g}
 g_{k\eta}(\xvec) = \frac{1}{(f(\xvec)^2 + \eta^2)^k}
\end{equation}
where $\eta$ is a real parameter and $k$ is integer. Since by hypotesis $\xvec_0$ is the only root of $f$, these functions have an absolute maximum at that point, where they reach the value $\frac{1}{\eta^{2k}}$.

We can now switch to the problem of finding the position of the maximum of $g_{k\eta}(\xvec)$, which is indeed equivalent to finding the zero of $f$. To solve this problem we borrow a very well known concept from classical mechanincs: the center of mass of a body with nonuniform mass density tends to lie near the points where the mass density is higher. With this in mind we can intuitively think that if a body has a mass density which diverges at a point, the center of mass should reasonably be near that point. Thus we regard the function $g_{k\eta}(\xvec)$ as the mass density of an $n$-dimensional object, whose shape is the domain $\Omega$, and compute the position of its center of mass when the parameter $\eta$ becomes infinitely small.

We recall that the center of mass $\xvec_{CM}$ of an $n$-dimensional non-uniform body is computed as
\begin{equation}
 \xvec_{CM} = \frac{\int_\Omega\xvec \rho(\xvec)\dnx}{\int_\Omega \rho(\xvec)\dnx}
\end{equation}
where $\rho(\xvec)$ is the mass density.
In our case we take $\rho(\xvec) = g_{k\eta}(\xvec)$ in the limit $\eta \rightarrow 0$, thus we write
\begin{equation}
\label{B}
 \xvec_{CM} = \mathlim{\eta}{0}\frac{\int_\Omega\frac{\xvec \dnx}{(f(\xvec)^2 + \eta^2)^k}}{\int_\Omega \frac{\dnx}{(f(\xvec)^2 + \eta^2)^k}}
\end{equation}
The key point is that the integrals in (\ref{B}) can always be made divergent in the limit $\eta \rightarrow 0$ for sufficiently high values of $k$, depending on the dimension of $\Omega$ and the multiplicity of the root.
Thus assuming a suitable value for $k$ we can restrict the integration to any domain over which the integrals diverge without affecting the result, in partucular we choose to integrate over $B^\varepsilon_{\xvec_0}$ and then perform a variable change $\xvec - \xvec_0 \rightarrow \xvec$, obtaining
\begin{equation}
\label{C}
\xvec_{CM} = \xvec_0 + \mathlim{\eta}{0}\frac{\int_{B^\varepsilon_0}\frac{\xvec \dnx}{(f(\xvec + \xvec_0)^2 + \eta^2)^k}}{\int_{B^\varepsilon_0} \frac{\dnx}{(f(\xvec + \xvec_0)^2 + \eta^2)^k}} = \xvec_0 + \mathlim{\eta}{0}\Delta \xvec_\eta
\end{equation}
Where $B_0^\varepsilon$ is the $\varepsilon$-set having the point $\xvec = 0$ on its boundary.
At this stage we can already see that the result of the computation is the searched solution $\xvec_0$ plus an $\eta$-dependent error term $\Delta \xvec_\eta$ which we shall show to go to zero as $\eta$ vanishes.

Since $B_0^\varepsilon$ is arbitrarily small we can subsitutue $f(\xvec + \xvec_0)$ with its lowest order Taylor expansion, which will depend on the multiplicity of the root. Without loss of generality we will assume that $f = |\xvec - \xvec_0|$ over $B_{\xvec_0}^\varepsilon$. In doing so we notice that, in the limit of vanishing $\eta$, at the denominator we are integrating a function with a singularity of the type $\frac{1}{|\xvec|^{2k}}$, while at the numerator the integrand has a $\frac{1}{|\xvec|^{2k-1}}$ singularity. The limit of the ratio of the two integrals is then zero and we can conclude that
\begin{equation}
\label{A}
 \xvec_{CM} = \mathlim{\eta}{0}\frac{\int_\Omega\xvec g_{n\eta}(\xvec)\dnx}{\int_\Omega g_{n\eta}(\xvec)\dnx} = \xvec_0
\end{equation}
Thus the position of the root of $f$ can be directly computed and corresponds to the center of mass of an $\Omega$-shaped object having a singular mass density.

\section{Numerical implementation}
Equation (\ref{A}) is exact and valid under broad assumptions and many numerical methods could be derived from it depending on how the integrals are actually handled. 
Any practical numerical implementation of (\ref{A}) will be confronted with the fact that the integrals are divergent in the limit $\eta \rightarrow 0$, a feature which is a necessary condition for the result to hold.
Even when using finite values of $\eta$ we are left with the problem of computing $n + 1$ $n$-dimensional integrals, a task whose computational cost could make the proposed method not competitive with other approaches, in fact, any deterministic discretization of the integrals would imply a computational cost growing very rapidly with the number of dimensions of $\Omega$. On the other hand statistical methods, such as Monte Carlo sampling, are much more effective for high dimensional problems because their rate of convergence is always $\frac{1}{\sqrt{N}}$ for $N$ samplings, regardless of the dimension\cite{MC1,MC2}.
Although a simple Monte Carlo sampling would then suffice, its rate of convergence would still be poor confronted with, e.g., the bisection method which has a rate of convergence of $\frac{1}{2^N}$ for $N$ iterations, or Newton-Raphson methods which are even faster, provided that the initial guess is close enough to the solution to make them converge.

\subsection{Parallel Monte Carlo sampling}

In computing the integrals in (\ref{A}) we can exploit a feature that will help speeding up things. We note that the two integrals to be computed are very similar, the only difference being that the function $g_{n\eta}$ is multiplied by $\xvec$ in the numerator. Moreover we know that the functions to be integrated reach very high values on a little region around the solution, while they are almost null everywhere else.
Both features can be taken into account by developing an adeguate Monte Carlo sampling procedure
in which the two integrals are computed in parallel: every time the numerator is sampled at the point $\xvec_i$, the denominator is sampled at the same point. 
Thus after $N$ samplings the estimated solution $\tilde\xvec_N$ is computed as
\begin{equation}
\label{MC}
 \tilde\xvec_N = \frac{\sum_{i=1}^N\frac{\xvec_ig_{k\eta}(\xvec_i)}{P(\xvec_i)}}{\sum_{i=1}^N\frac{g_{k\eta}(\xvec_i)}{P(\xvec_i)}}
\end{equation}
where $\xvec_i$ is a sequence of points in $\Omega$, generated according to the probability density $P$, which is used for the evaluation of both the numerator and the denominator.
Thus if during the sampling the numerator has a sudden jump due to $\xvec_i$ lying near the singularity, the denominator will jump too, and this will reduce the fluctuations of the ratio, causing $\tilde\xvec_N$ to jump towards the exact value. This parallel-sampling choice has no formal justification besides the fact that it turns out to work well, in fact, even if we know that both integrals will converge to the respective solutions as $\frac{1}{\sqrt{N}}$, we cannot in principle tell anything about the rate of convergence of their ratio, because the two samplings are higly correlated in way that depends on the particular problem.
Another useful property of (\ref{MC}) derives from the fact that we are computing the ratio of two integrals and thus the final result is independent from the normalization of $P(\xvec)$, which is a major issue in Monte Carlo sampling techniques.

\subsection{Uniform sampling}
The simplest choice is to adopt a constant probability density $P(\xvec_i) = 1$ to generate the sampling points. 
In the following we shall show that in this case the method converges even in the case $\eta = 0$, and its rate of convergence is $\frac{1}{N^{\frac{1}{n}}}$.
To do so we evaluate (\ref{C}) using the parallel-sampling method with uniform probability density, after setting $\eta = 0$. We obtain
\begin{equation}
 \tilde\xvec_N = \xvec_0 + \frac{\sum_{i=1}^N\frac{\xvec_i}{f(\xvec + \xvec_0)^{2k}}}{\sum_{i=1}^N\frac{1}{f(\xvec + \xvec_0)^{2k}}} = \xvec_0 + \Delta \xvec_N
\end{equation}
where $\xvec_i$ is a sequence of uniformly distributed points inside $B^\varepsilon_0$ and $\Delta \xvec_N$ is the error after $N$ samplings.

If we now consider the quantity $|\Delta \xvec_N|$ we can write
\begin{equation}
\label{D}
 |\Delta \xvec_N| \leq \frac{\sum_{i=1}^N\frac{1}{|\xvec_i|^{m-1}}}{\sum_{i=1}^N\frac{1}{|\xvec_i|^m}}
\end{equation}
since $f(\xvec + \xvec_0)^{2k} \propto |\xvec|^m$, for some $m$, on an arbitraryly small domain around $\xvec = 0$.

We can rearrange the elements of the sequence $|\xvec_i|$ in order to $|\xvec_1|$ to be the smallest value and then write 
\begin{equation}
\label{E}
 \sum_{i=1}^N\frac{1}{|\xvec_i|^m} = \frac{1}{|\xvec_1|^m} + \sum_{i=2}^N\frac{1}{|\xvec_i|^m}
\end{equation}
Now substituting (\ref{E}) in (\ref{D}) we obtain
\begin{equation}
 |\Delta \xvec_N| \leq \frac{\frac{1}{|\xvec_1|^{m-1}} + \sum_{i=2}^N\frac{1}{|\xvec_i|^{m-1}}}
{\frac{1}{|\xvec_1|^m} + \sum_{i=2}^N\frac{1}{|\xvec_i|^m}}\;. 
\end{equation}
Since the points $\xvec_i$ are uniformly distributed over the $n$-dimensional domain $\Omega$ the average distance between two points, in the limit $N \rightarrow \infty$, will be of the order of $\frac{1}{N^\frac{1}{n}}$ and thus $|\xvec_1|$ will be of that order too. Then, in the limit of large $N$, the term $1/|\xvec_1|^m$ in (\ref{E}) will diverge while the rest of the summation will remain finite, we can then say that
\begin{equation}
 \mathlim{N}{\infty}|\Delta \xvec_N| 
\leq \mathlim{N}{\infty}\frac{\frac{1}{N^\frac{m-1}{n}}}{\frac{1}{N^\frac{m}{n}}} 
= \mathlim{N}{\infty}\frac{1}{N^\frac{1}{n}} = 0
\end{equation}
and thus the sequence $\tilde\xvec_N$ converges to $\xvec_0$ with a rate of convergence of $\frac{1}{N^\frac{1}{n}}$.

Thus in the case of uniform sampling the method is (almost surely) convergent but its rate of convergence is poor, depends on the dimension of the space and in particular becomes lower for higher dimensional problems.

\subsection{Adaptive sampling}
The efficiency of Monte Carlo sampling techniques is highly enhanced when the distribution of sampling points resembles the function to be integrated, thus we developed an adaptive-sampling method in which the probability density is varied during the process in order to maximize the efficiency.
The key feature that allows such a strategy is again the fact that we are computing the ratio of two integrals and thus we are allowed to use non-normalized probability densities for the generation of sampling points. In particular, the parallel-sampling tecnique allows us to change the probability during the process: a peaked probability distribution is chosen for the generation of the sampling points, and its shape is gradually varied in order to concentrate the sampling around the solution. 

The method proceeds as follows: a gaussian probability density of the form
\begin{equation}
 P^i(x_1,...,x_n) = \prod_{j=1}^ne^\frac{(x_j - \tilde x_j^i)^2}{(\sigma_j^i)^2}
\end{equation}
is used for the distribution of sampling points, where $x_j$ is the $j$th-component of $\xvec$, $\tilde x_j^i$ is the estimation of the solution at step $i$ and $\sigma_j^i$ is a suitable measure of the fluctuations of $\tilde x_j$ up to step $i$.
When the process starts, $\sigma_j^0$ is set to a value such that the gaussian is flat over $\Omega$.
The estimation $\tilde \xvec^0$ loses all meaning in this case and can be set to an arbitrary point.
After a few samplings, a new estimation of the solution is available, as well as an estimation of its fluctuations. These values are then used to update the probability density $P^i$. The key point is that as the sampling proceeds, the estimation comes closer to the exact solution and the fluctuations are reduced, thus $P^i$ is gradually sharpened and centered around $\xvec_0$. This makes $P^i$ resemble the functions to be integrated, making the sampling more efficient and further reducing the fluctuations, which in turn make $P^i$ more peaked around the solution.

This adaptive-sampling technique must be tuned carefully. If the variance of the gaussian is reduced too fast, it can happen that the region around $\xvec_0$ is never sampled and thus the process converges to a random point and fails to find the correct solution. Conversely, if the variance is reduced too slowly, ther is no real advantage in using an adaptive algorithm and the results are the same as using a constant uniform probability density.
In practice the key parameters are the number of samplings between two successive updates of $P^i$ and the way $\sigma$ is computed, a suitable choice has been found to be
\begin{equation}
\label{sigma}
\sigma_j^i = \frac{1}{N}\sum_{k=0}^{N-1}|\tilde x_j^{i-k} - \tilde x_j^{i-k-1}| 
\end{equation}
where the number of samplings $N$ over wich the mean is performed determines the speed at which the variance of the gaussian is reduced, which in turn affects the rate of convergence.

\begin{figure}
 \includegraphics*[width=\linewidth,height=\linewidth]{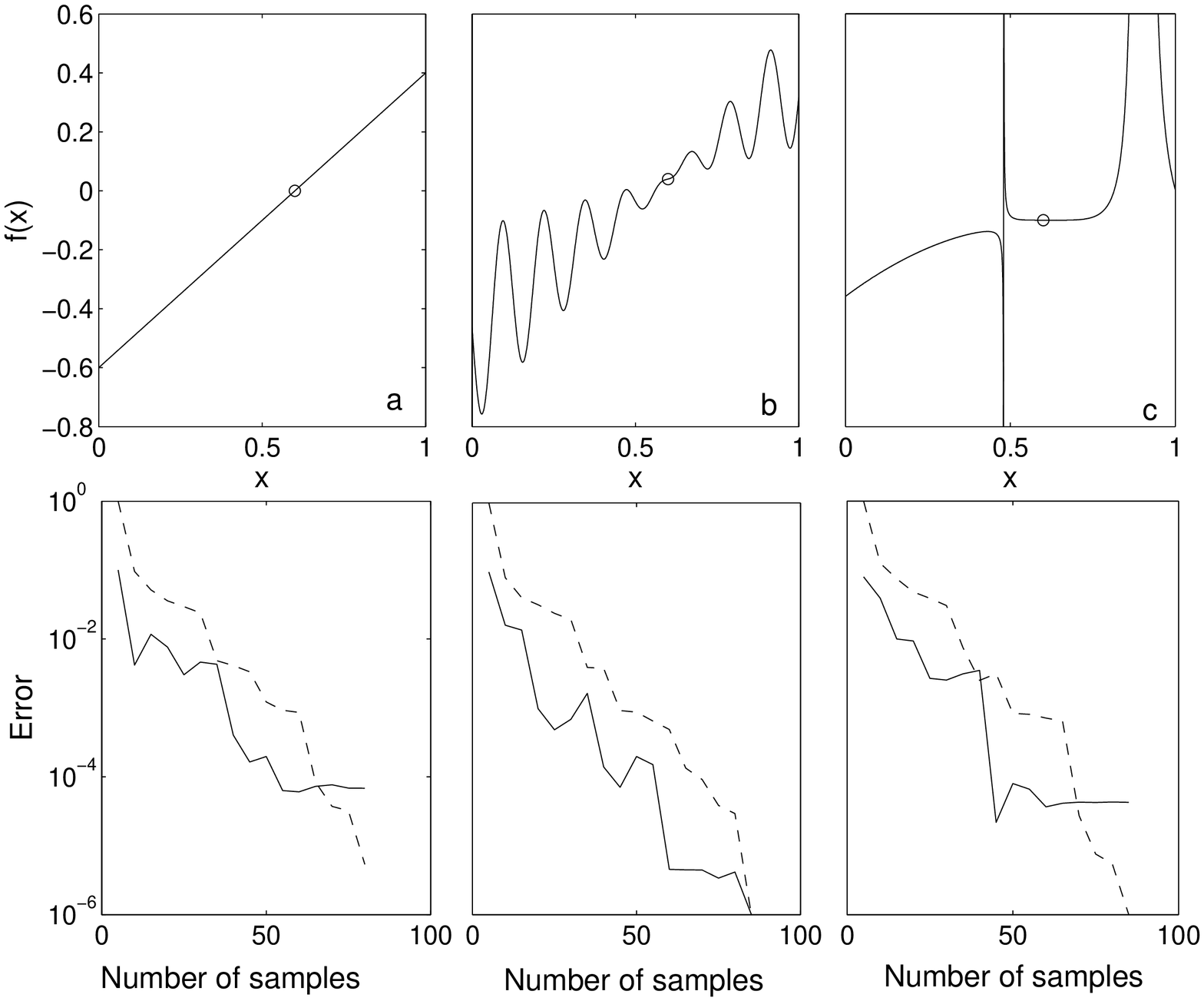}
\caption{Three one dimensional root-finding test problems solved by the adaptive Monte Carlo sampling method. Figures on the top show the test functions, each of which has a zero at $x_0 = 0.6$. Figures on the bottom row show the corresponding convergence tests for the adaptive Monte Carlo sampling method: dashed lines show the the fluctuation $\sigma$ of the solution, the continuous line is the actual error $|\tilde x - x_0|$. }\label{fig1}
\end{figure}

Figure (\ref{fig1}) shows three one-dimensional test problems (top figures) against which the method has been tested in order to check its convergence.
In these examples $P$ is updated every $N = 5$ samplings and $\eta = 10^{-8}$
The graphs at the bottom show typical convergence curves from which it can be seen that the rate is approximatively exponential (continuous line).
It can also be noted that the parameter $\sigma$ (dashed line) usually overestimates the actual error and can thus be taken as a conservative estimate.
Actually, Figs. \ref{fig1}a and \ref{fig1}c show exceptions to this behaviour: this is caused by the finite value of $\eta$ which introduces an intrinsic uncertainty in the determination of the solution and thus the error curves level at a certain value; anyway such error can be made arbitrarily small by reducing $\eta$.
\begin{figure}
  \includegraphics*[width=\linewidth,height=\linewidth]{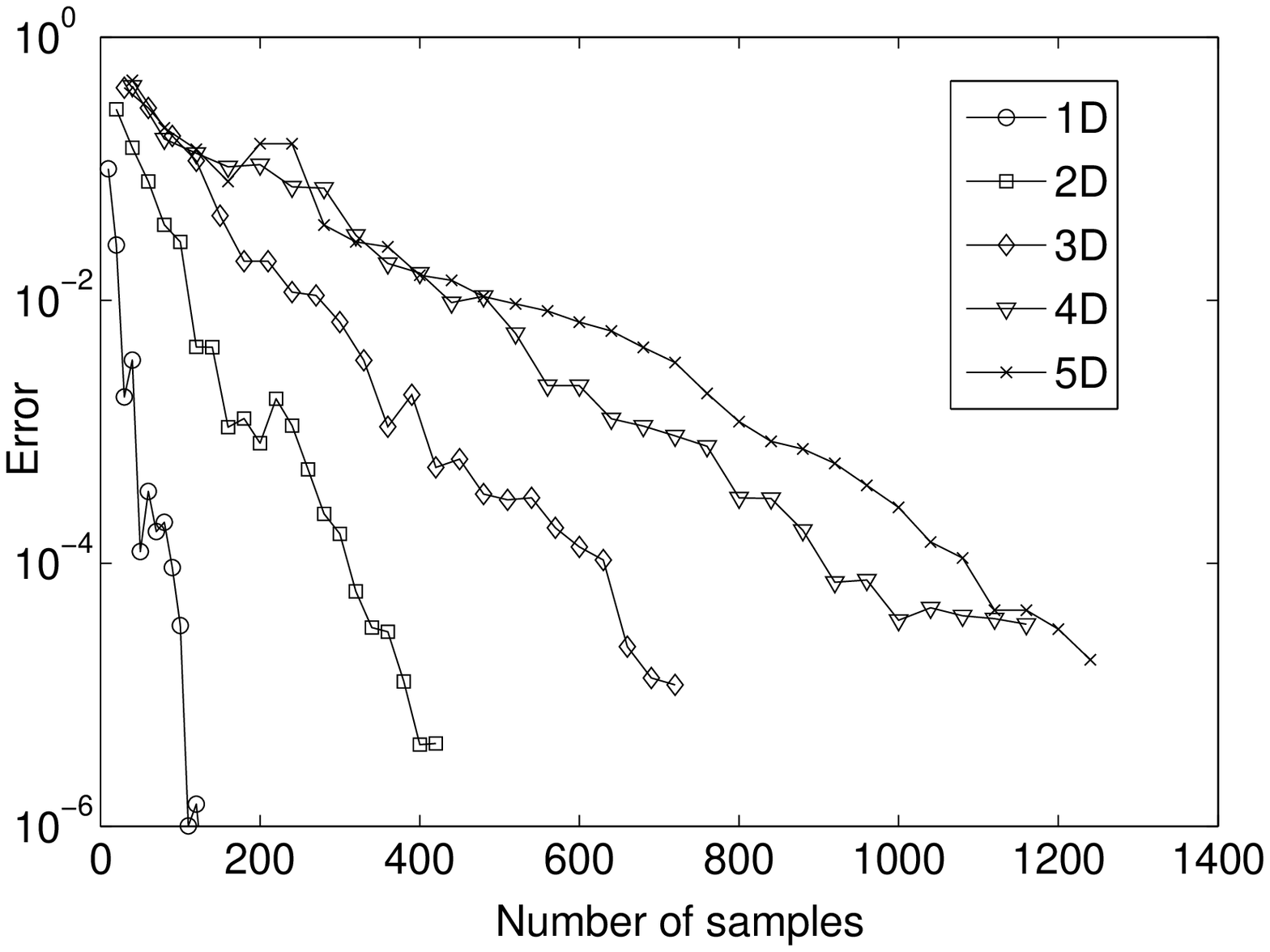}
\caption{Convergence test for multidimensional problems. The rate of convergence is lower for high dimensional cases, but remains exponential. }\label{fig2}
\end{figure}
A similar test has been performed for multidimensional problems in Fig. \ref{fig2} where the error $|\tilde \xvec - \xvec_0|$ is plotted against the number of Monte Carlo samplings. The function considered is the $n$-dimensional equivalent of Fig. \ref{fig1}b.
It can be seen that the exponential rate of convergence is mantained in every dimension, differing only in the slope of the function. This is due to the fact that as the number of dimensions is increased, the parameter $N$ in (\ref{sigma}) has to be increased in order to ensure convergence. In this case $N = 10$ for the 1D problem and $N = 50$ for the 5D problem.

\subsection{Automated root search}
We have shown that the adaptive sampling method dramatically increases the rate of convergence with respect to the uniform sampling, at the cost of introducing the possibility that it could fail to converge in some cases.
This weakness can actually be effectively turned into a benefit.
The failure of the adaptive-sampling tecnique derives from having introduced some locality in the method, indeed, using a gaussian probability density corresponds to limiting the domain considered for the computation of integrals in (\ref{A}).
This locality can then drive the method to converge to local minima or to points inside regions where $f$ reaches very low values.
Actually the worst case would be that of a function $f$ which has two roots in $\Omega$: in this case the function $g_{k\eta}$ would have two singularities that would make the uniform-sampling method not converge to any of the two.
On the other side, the adaptive-sampling method, due to its tendency to drive the probability density around points where $f$ has very low values, will converge to one of the two solutions.

This feature can be exploited to develop an automated algorithm that sequentially discovers all the roots of $f$. 
Consider a function $f$ defined on $\Omega$ which has a number of roots at points $\xvec_1, ... , \xvec_M$.
The adaptive-sampling method then the method will converge to one of the solutions, say $\xvec_1$.
A this point if the region surrounding $\xvec_1$ is removed from the domain $\Omega$ and the process is repeated, the method will converge to another solution, say $\xvec_2$. This loop can be iterated in order to find all the roots of $f$.

\section{Summary and conclusions}
We have proposed a global root-finding algorithm which allows to compute the position of the root of a scalar function over an $n$-dimensional domain under very broad assumptions.
An exact analytical result is given which transforms the inverse problem of finding the root of a function to a direct computation of its position, interpreted as the center of mass of an $n$-dimensional $\Omega$-shaped object with singular mass density.
In view of this result a parallel Monte Carlo sampling method has been developed and proven to converge with a rate of $\frac{1}{N^\frac{1}{n}}$.
The method is then extended in order to employ adaptive sampling of the integrals, providing a numerical method with exponential rate of convergence which can be further extended to an automated root-search algoritm applicable to functions with multiple roots.

The proposed methods, being based on integral quantities, are global in nature and applicable to a wide variety of cases in which the function is bad behaved (not differentiable, singular, rapidly oscillating, ...) or the domain $\Omega$ is not regular (not simply connected, irregular shape, ...).

This paper is a demonstration of the possibilities of the proposed method, further work is in progress in order to improve its convergence rate and range of applicability.
A wide variety of methods can be derived from equation (\ref{A}), depending on the specific choice of the $g_{k\eta}$ function and the specific method used to compute the integrals, and this could open the way towards more efficient methods.

We are grateful to Fausto Rossi and David Taj for stimulating and fruitful discussion.

\end{document}